\documentclass[11pt,a4paper]{article}
\usepackage[russian,english]{babel}
\usepackage{amssymb}
\begin{document}
\title{Global Structure of Locally Convex Hypersurfaces in Finsler-Hadamard Manifolds}
\author{A. A. Borisenko, E.A. Olin}
\date{21 March 2008}

\maketitle

\textbf{Abstract}

Locally convex compact immersed hypersurfaces in Finsler-Hadamard manifolds with bounded $\mathbf{T}$-curvature are considered. We prove that under certain conditions on the normal curvatures such hypersurfaces are embedded as the boundary of convex body and are homeomorphic to the sphere.
\section{Introduction}
Let $M$ be a complete Finsler manifold. Then
\begin{enumerate}
 \item
A set $A$ is said to be  \textit{convex}, if each shortest path with endpoints in $A$ is entirely contained in $A$.

 \item A set $A$
is said to be \textit{locally convex}, if each point $P \in A$ has a neighborhood $U_P$ in $M$ such that the set $A\cap
U_P$ is convex.
\end{enumerate}

Hadamard proved the following theorem.

\textbf{Theorem. [10]} \textit{Let $\varphi$ be an immersion of compact $n$-dimensional oriented manifold $M$ in Euclidean space $E^{n+1}$, $n\geqslant2$ with everywhere positive Gaussian curvature. Then $\varphi(M)$ is a convex hypersurface.}

Chern and Lashof generalized this theorem.

\textbf{Theorem. [10]} \textit{Let $\varphi$ be an immersion of compact $n$-dimensional oriented manifold $M$ in Euclidean space $E^{n+1}$, $n\geqslant2$. Then the following two assertions are equivalent.}
\begin{enumerate}
\item \textit{The degree of the spherical mapping equals $\pm 1$ and the Gaussian curvature does not change its sign;}

\item \textit{$\varphi(M)$ is a convex hypersurface.}
\end{enumerate}

A topological immersion $f:N^n\rightarrow M^{n+1}$ of a manifold
$N^n$ into a Riemannian manifold $M^{n+1}$ is called \textit{locally convex} at a point  $x\in N^n$ if $x$ has a neighborhood $U$
such that $f(U)$ is a part of the boundary of a convex set in
$M^{n+1}$.

Heijenoort proved the following theorem.

\textbf{Theorem. [9] }\textit{Let $f:N^n\rightarrow E^{n+1}$,
$n\geqslant 2$ be a topological immersion of a connected manifold $N^n$. If $f$ is locally convex at all points and has at least one point of local strict support and $N^n$
is complete in the metric induced by this immersion, then $f$ is an embedding and $f(N^n)$ is the boundary of a convex body.}

S. Alexander [1] (see also A. A. Borisenko [5]) generalized Hadamars's theorem for compact immersions when an ambient space is a complete simply connected manifold of non-positive curvature (\textit{Hadamard} manifold).

\textbf{Theorem. [1,5] }\textit{Let $f:N^n\rightarrow M^{n+1}$,
$n\geqslant 2$ be an immersion of compact connected manifold
$N^n$ in a complete simply connected Riemannian manifold $M^{n+1}$ of non-positive sectional curvature. If the immersion $f$ is locally  convex, then $f$ is an embedding, $f(N^n)$ is the boundary of a convex set in  $M^{n+1}$, and $f(N^n)$ is homeomorphic to the sphere
$\mathbb{S}^n$.}

The goal of this paper is to generalize this theorem to an immersion of compact manifold into a complete simply connected Finsler manifold of non-positive curvature (\textit{Finsler-Hadamard} manifold).

\textbf{Theorem 1.} \textit{Let $f:N^n\rightarrow M^{n+1}$, $n\geqslant 2$ be an immersion of a compact connected
manifold $N^n$ in a complete simply connected Finsler manifold $M^{n+1}$. Let $N^n$ and $M^{n+1}$ satisfy the following conditions:
 \begin{enumerate}
 \item The flag curvature $K \leqslant -k^2$, $k\geqslant 0$;
 \item The $\mathbf{T}$-curvature $|\mathbf{T}| \leqslant \delta$, where $0\leqslant\delta < k$;
 \item All the normal curvatures of
$N^n$  $\mathbf{k}_{\mathbf{n}} > 2 \delta$.
\end{enumerate}
Then $f$ is an embedding, $f(N^n)$ is the
boundary of a convex set in  $M^{n+1}$ which is homeomorphic to the ball, and $f(N^n)$ is homeomorphic to the sphere $\mathbb{S}^n$.}

We also show that the theorem of S. Alexander holds for Berwald spaces without any additional restrictions.

Note, that if all the flag curvatures vanish (the case when $k=0$) then $\mathbf{T} = 0$ and we obtain the generalization of Hadamard's theorem to Minkowski spaces.

At the end of the paper we give an example which is based on the recent work of V. K. Ionin [13]. This example shows the necessity of the ambient space to be complete.

\section{Preliminaries}

In this section we recall some basic facts and theorems from
Finsler geometry that we need. See [11] for details.

 Let $M^n$ be an
$n$-dimensional connected $C^{\infty}$-manifold. Denote by
$TM^n=\bigsqcup_{x\in M^n}T_xM^n$ the tangent bundle of $M^n$,
where $T_xM^n$ is the tangent space at $x$. A \textit{Finsler
metric} on $M^n$ is a function $F:TM^n\rightarrow [0,\infty)$ with
the following properties:
\begin{enumerate}
\item $F\in C^{\infty}(TM^n\backslash\{0\})$;

\item $F$ is positively homogeneous of degree one, i. e. for any
pair $(x,y)\in TM^n$ and any $\lambda>0$, $F(x,\lambda y)=\lambda
F(x,y)$; \item For any pair $(x,y)\in TM^n$ the following bilinear
symmetric form $g_y:T_x M^n\times T_x M^n\rightarrow \mathbb{R}$
is positively definite,
$$\mathbf{g}_y(u,v):=\frac{1}{2}\frac{\partial^2}{\partial t \partial s}
\lbrack F^2(x,y+su+tv)\rbrack |_{s=t=0}$$
\end{enumerate}

The pair $(M^n,F)$ is called \textit{a Finsler manifold}.

If we denote by $$\mathbf{g}_{ij}(x,y) =
\frac{1}{2}\frac{\partial^2}{\partial y^i\partial
y^j}[F^2(x,y)],$$ then one can rewrite the form
$\mathbf{g}_y(u,v)$ as
$$\mathbf{g}_y(u,v)=\mathbf{g}_{ij}(x,y)u^iv^j$$

For any fixed vector field $Y$ defined on the subset $U \subset M^n$, $\mathbf{g}_Y(u,v)$ is a Riemannian
metric on $U$.

Given a Finsler metric $F$ on a manifold $M^n$. For a smooth curve
$c:[a,b]\rightarrow M^n$ the length is defined by the integral
$$L_F(c)=\int_a^b F(c(t),\dot{c}(t))dt = \int_a^b \sqrt{\mathbf{g}_{\dot{c}(t)}(\dot{c}(t),\dot{c}(t))}dt.$$

Here the last integral is over $c(t)$.

Chern covariant derivative
$$\nabla:T_xM^n \times C^\infty(TM^n)\rightarrow T_xM^n$$
is defined as follows
$$\nabla_yU:=\{dU^i(y)+U^jN_j^i(x,y)\}\frac{\partial}{\partial x^i}|_x,$$
where $y\in T_xM^n$,  $U\in C^\infty(TM^n)$ è
$$N^i_j(x,y)=\frac{\partial}{\partial y^j}\left[\frac{1}{4}
g^{il}(x,y)\left\{2\frac{\partial g_{ml}}{\partial
x^k}(x,y)-\frac{\partial g_{mk}}{\partial x^l}(x,y)
\right\}y^my^k\right].$$

Defined connection is not an affine connection for general Finsler metrics. The connection will be affine if and only if (see [11, p.86]) given Finsler metric is a \textit{Berwald metric}. By definition this means that the geodesic equations have the same form as in Riemannian geometry, or the geodesic coefficients

$$G^i(x,y)=\frac{1}{4}
g^{il}(x,y)\left\{2\frac{\partial g_{jl}}{\partial
x^k}(x,y)-\frac{\partial g_{jk}}{\partial x^l}(x,y)
\right\}y^jy^k$$ can be expressed as
$$G^i(x,y) = \Gamma^i_{jk}(x)y^j y^k. $$

 For a vector $y\in T_xM^n
\backslash \{0\}$ consider the functions
$$R_k^i(y)=2\frac{\partial G^i}{\partial x^k}-\frac{\partial^2 G^i}{\partial x^j \partial
y^k}y^j+2G^j\frac{\partial^2G^i}{\partial y^j \partial
y^k}-\frac{\partial G^i}{\partial y^j}\frac{\partial G^j}{\partial
y^k}.$$
The family of linear transformations
$$\mathbf{R}=\left\{\mathbf{R}_y=R_k^i(y)\frac{\partial}{\partial x^i}\otimes dx^k|_x:T_xM^n\rightarrow T_xM^n,y\in T_xM^n\backslash\{0\}, x \in
M^n\right\}$$ is called the \textit{Riemannian curvature} [11, p.97].

Let $P\subset T_xM^n$ be a tangent 2-plane. For a vector $y \in P\backslash\{0\}$ define
$$K(P,y):=\frac{\mathbf{g}_y(\mathbf{R}_y(u),u)}{\mathbf{g}_y(y,y)\mathbf{g}_y(u,u)-\mathbf{g}_y(y,u)^2},$$
where $P=\mathrm{span}\{y,u\}$. $K(P,y)$
does not depend on $u\in P$, but depend on $y$. The value $K(P,y)$ is called the
\textit{flag curvature} of the flag $(P,y)$ in $T_xM^n$.

 Then a complete simply connected Finsler manifold of non-positive curvature is called(\textit{Finsler-Hadamard} manifold).
In these spaces Cartan-Hadamard theorem holds [8].

Consider a geodesic $c:[a,b]\rightarrow M^n$. Then a
\textit{Jacobi field} along the geodesic $c$ is a vector filed
$J(t)$ that satisfies the \textit{Jacobi equation}:
$$ \nabla_{\dot{c}(t)}\nabla_{\dot{c}(t)}J+\mathbf{R}_{\dot{c}(t)}(J)=0 $$

 The flag curvature does not describe all the properties of Finsler spaces. Therefore so-called none-Riemannian curvatures are considered. They are equal to 0 for Riemannian spaces. We will need one of them which is called the
$\mathbf{T}$-curvature [11, p.153].

Let $(M^n,F)$ be a Finsler space. For a given vector $y \in
T_xM^n\backslash\{0\}$ denote by $Y$ its extension to  a geodesic field in a
neighborhood of $x$. Let $\nabla$ denote the Chern connection,
$\tilde{\nabla}$ denote the Levi-Civita connection of the induced
Riemannian metric $\tilde{g}=\mathbf{g}_Y$. For a vector $v \in
T_xM^n$ define
\begin{equation}
\mathbf{T}_y(v)=\mathbf{g}_y(\nabla_v
V,y)-\tilde{g}(\tilde{\nabla}_v V,y)
\end{equation}
 where $V$ is a vector field such that $V_x=v$.

The function $\mathbf{T}_y(v)$, $y \in T_xM^n\backslash\{0\}$ is
called $\mathbf{T}$\textit{-curvature}.

$\mathbf{T}$-curvature is said to be bounded above
$\mathbf{T}\geqslant-\delta$ if [11, p.223]

$$\mathbf{T}_y(u)\geqslant-\delta\left[\mathbf{g}_y(u,u)-\mathbf{g}_y\left(u,\frac{y}{F(y)}\right)^2\right]F(y)$$

The upper bound is defined at the same manner.

Notice that the  $\mathbf{T}$-curvature
vanish for Berwald metrics; the converse is also true [11, p.155].

Let $\varphi :N\rightarrow M^n$ be a hypersurface in
$M^n$. A vector $\mathbf{n} \in T_{\varphi(x)}M^n$ is called normal vector to $N$ at the point $x \in N$ if
$\mathbf{g}_{\mathbf{n}}(y,\mathbf{n})=0$ for each $y \in T_xN$.
It is known that such a vector exists [11, p.27]. Note that for non-reversible metrics the vector $-\mathbf{n}$ is not the normal vector.

We can consider the subbundle $\nu(N)$ of the tangent bundle $TM^{n}$ formed by all normal to $N$ vectors with chosen orientation; $\nu(N)$ is called the \textit{normal bundle} over $N$.

The map
 $\exp_{N}:\nu(N)\rightarrow M^n$
defined by
$$\exp_{N}(x,\mathbf{n})=\exp_x(\mathbf{n}) $$
ia called an exponential map of a hypersurface $N$.

For a hypersurface $N$ in $M^{n}$ the
\textit{normal curvature} $\mathbf{k}_\mathbf{n}$ at a point $x\in N$ in a direction $y\in T_xN$ is defined by
$$\mathbf{k}_\mathbf{n} = \mathbf{g}_\mathbf{n}(\nabla_{\dot{c}(t)}\dot{c}(t)|_{t=0},\mathbf{n}),$$
where $c(t)$ is a geodesic in the induced connection on $N$ such that $\dot{c}(0)=y$ and $\mathbf{n}$ is unit normal vector.

\textbf{Proposition 1. [11, p.217]}\textit{ If at a point $x \in N$ all the normal curvatures $\mathbf{k}_\mathbf{n} > 0$ or $\mathbf{k}_\mathbf{n} < 0$ then $N$ is locally convex at $x$. }

As  we have noticed, for any fixed vector field $Y$ defined on the subset $U \subset M^n$, $\mathbf{g}_Y(u,v)$ is a Riemannian
metric on $U$. This metric satisfies a number of useful properties.

\textbf{Proposition 2. [11, p.100]}\textit{ For a non vanishing geodesic vector field $Y$,
consider the induced Riemannian metric $\tilde{g} =
\mathbf{g}_{Y}$. Define by $K(P,Y)$ the flag curvature of the initial Finsler metric, and by $\tilde{K}(P,Y)$ the sectional curvature of the induced metric. Then for each flag $P$: }
$$K(P,Y) = \tilde{K}(P,Y).$$

For a function $\rho$ on $M^n$ at a point $x\in T_x M^n$ the gradient $\mathbf{grad} (\rho)_x \in T_x^*M^n$ is defined as follows [11, p.41]: $$\frac{\partial \rho}{\partial x^i}(x)v^i=\mathbf{g}_{\mathbf{grad} (\rho)_x}(\mathbf{grad} (\rho)_x,v), v\in T_xM^n.$$

\textbf{Proposition 3. [11, p.216]}\textit{ Let $\rho$ be a $C^{\infty}$ distance function defined on an open subset $U\subset
M^n$. Consider the induced Riemannian metric $\tilde{g} =
\mathbf{g}_{\mathbf{grad} (\rho)}$. Let $\mathbf{k}_{\mathbf{n}}$ and
$\tilde{\mathbf{k}}_{\mathbf{n}}$ respectively be the normal curvatures of $N=\rho^{-1}(s)$ in the metrics $F$ and $\tilde{g}$
with respect to the unit normal vector $\mathbf{n}=\mathbf{grad}(\rho)_x$, $x\in N$. Then for a vector
$y\in T_xN$:}

$$\mathbf{k}_{\mathbf{n}}(y)=\tilde{\mathbf{k}}_{\mathbf{n}}(y) -
\mathbf{T}_{\mathbf{n}}(y)$$

\textbf{Proposition 4. [11, p.235]} \textit{Let $(M^n,F)$  be a complete simply connected Finsler manifold with the flag
curvature $K\leqslant-k^2$ and the  $\mathbf{T}$-curvature $|\mathbf{T}|\leqslant\delta$ such that  $\delta > k$. Then
the balls of arbitrary radii in $M^n$ are convex.}

For a hypersurface $N$ we can define the \textit{shape operator} as follows. Let $\rho$ be a $C^{\infty}$ distance function from the hypersurface $N$ defined on a tubular neighborhood of $N$. Let $\tilde{\nabla}$ be the Levi-Civita connection of the induced Riemannian metric
$\tilde{g}=\mathbf{g}_{\mathbf{grad} (\rho)}$.  Then for the normal vector
$\mathbf{n}=\mathbf{grad}(\rho)|_N$ the operator $S_{\mathbf{n}}(w):T_xN\rightarrow T_xN$ defined by
$$S_{\mathbf{n}}(w)=\tilde{\nabla}_w \mathbf{n}$$
is called the\textit{ shape operator} [11, p.221].

The shape operator $S$ satisfies the following properties [11, p.222].

\begin{enumerate}
\item

 $\mathbf{g}_{\mathbf{n}}(S_{\mathbf{n}}(y),v) = \mathbf{g}_{\mathbf{n}}(y,S_{\mathbf{n}}(v))$

\item

 $\tilde{\mathbf{k}}_{\mathbf{n}}(y) =
\mathbf{g}_{\mathbf{n}}(S_{\mathbf{n}}(y),y)$

\item For the family of the shape operators $S_t$ of the hypersurfaces $N_t=\rho^{-1}(t)$ the Ricatty equation holds

$$\dot{S}_t+S_t^2+\mathbf{R}_{\dot{c}(t)}=0$$
Here $c(t)$ is the integral curve of $\mathbf{grad} (\rho)$.

\end{enumerate}

And the differential equation on the normal curvatures of the parallel hypersurfaces can be obtained in the form

$$\dot{\tilde{\mathbf{k}}}+ \tilde{\mathbf{k}}^2 + f(t) = 0$$

\section{Comparison theorem for Jacobi fields}
The Rauch comparison theorem was proved for comparison of the lengths of Jacobi fields along geodesics in different Riemannian manifolds. Berger (for example see [6]) has given an extension of Rauch's theorem for the exponential map of the geodesics. Warner [12] proved such a result for Jacobi fields associated with submanifolds.

Here we extend Rauch's theorem for Jacobi fields associated with hypersurfaces in Finsler manifold.

 Consider the smooth hypersurface $N$ in $M^n$. Consider
the geodesic $c(t)$ started at a point $x \in N$ in the normal direction to $N$ at $x$.

We will call the Jacobi fields $J(t)$ along the geodesic $c(t)$ a
$N$\textit{-Jacobi field} if
$$\nabla_{\dot{c}(t)}J(t)|_{t=0} = S_{\dot{c}}(J(t))|_{t=0}$$

A point $c(t_0)$ is called a focal point to $N$ along $c$ if
there exist non-trivial $N$-Jacobi field $J$ along $c$ such that
$J(t_0)=0$

As in Riemannian geometry focal point can be considered as critical value of the exponential map of the hypersurface.

Consider the index form with respect to a hypersurface $N$:
$$I_t(X,Y)=\mathbf{g}_{\mathbf{n}}(S_{\mathbf{n}}(X),Y)|_0+\int_0^t\left[\mathbf{g}_{\dot{c}}(\nabla_{\dot{c}}X,\nabla_{\dot{c}}Y)-\mathbf{g}_{\dot{c}}(\mathbf{R}_{\dot{c}}(X),Y)\right]ds$$

Note that the index form is symmetric, i.e.

$$I_t(X,Y)=I_t(Y,X)$$

The index form naturally arises from the second variation formula [11, p. 160].

We prove the extremal property of Jacobi fields (Index lemma).

\textbf{Lemma 1.}\textit{ Let $N$ be a smooth hypersurface in
Finsler manifold $M^n$, $c:[0,s]\rightarrow M^n$ be a normal geodesic such that $c(0) \in N$, $\dot{c}(0)$ be
$\mathbf{g}_{\dot{c}(0)}$-orthogonal to $N$ and there are no point focal to $N$ along $c$. Denote by $J$, $Y$ non-trivial vector fields such that $J(s) = Y(s)$ and field $J$ is $N$-Jacobi. Then}

$$I_s(J,J)\leqslant I_s(Y,Y)$$
and equality occurs if and only if $J=Y$.

\textit{Proof.} We follow the idea of [6] in the proof. Let
$J_1$,...,$J_{n-1}$ be a basis of the space of $N$-Jacobi fields along
$c$. Then we can express $J$ and $Y$ as follows $J(t) = a^iJ_i(t)$,
$Y(t) = y^i(t)J_i(t)$, where $a^i$ are constants and $t\neq 0$.

Denote $A(t)=y'^i(t)J_i(t)$ and $B(t) = y^i(t)J_i'(t)$. Then $Y' = A
+ B$ and we get

$$I_s(Y,Y) =\mathbf{g}_{\mathbf{n}}(S_{\mathbf{n}}(Y),Y) + \int_0^s\left[ \mathbf{g}_{\dot{c}}(A,A)+2\mathbf{g}_{\dot{c}}(A,B)+\mathbf{g}_{\dot{c}}(B,B)-\mathbf{g}_{\dot{c}}(\mathbf{R}_{\dot{c}}Y,Y)\right]dt$$

From the Lagrange equality [2, p.135] we obtain

$$\mathbf{g}_{\dot{c}}(J_i,J_j') - \mathbf{g}_{\dot{c}}(J_j,J_i') =const$$

The $N$-Jacobi condition implies

$$\mathbf{g}_{\dot{c}}(J_i,J_j') - \mathbf{g}_{\dot{c}}(J_j,J_i')
=0$$

Then

$$\mathbf{g}_{\dot{c}}(Y,B) ' = y'^i y^j\mathbf{g}_{\dot{c}}(J_iJ_j')+y^iy'^j\mathbf{g}_{\dot{c}}(J_i',J_j)+y^iy^j\mathbf{g}_{\dot{c}}(J_i',J_j')+y^iy^j\mathbf{g}_{\dot{c}}(J_i,J_j'').$$

and

$$I_s(Y,Y)=\mathbf{g}_{\mathbf{n}}(S_{\mathbf{n}}(Y),Y)|_{0} +
\mathbf{g}_{\dot{c}}(Y,B)|_0^s+\int_0^s
\mathbf{g}_{\dot{c}}(A,A)dt-\int_0^s
y^iy^j\mathbf{g}_{\dot{c}}(J_i,J_j''+\mathbf{R}_{\dot{c}}J_j)dt
$$

Using the $N$-Jacobi condition we get
$$\mathbf{g}_{\mathbf{n}}(S_{\mathbf{n}}(Y),Y)|_{0} -
\mathbf{g}_{\dot{c}}(Y,B)|_0 = 0$$

Next, since $J(s) = Y(s)$ then $a^i = y^i(s)$. And we have

$$\mathbf{g}_{\dot{c}}(Y,B)|_s = \mathbf{g}_{\dot{c}}(J,a^iJ_i')|_s =
\mathbf{g}_{\dot{c}}(J,J')|_s = I_s(J,J).$$

The reader should see [2, p.177] for the proof of last equality.
Hence,

$$I_s(Y,Y) = I_s(J,J) + \int_0^s
\mathbf{g}_{\dot{c}}(A,A)dt$$

And the lemma follows. $\square$

We consider two Finsler manifolds $M^n$ and $\bar{M}^n$. In the manifold
$M^n$ we consider a hypersurface  $N$, a normal geodesic  $c:[0,s]\rightarrow M^n$ such that $c(0) \in N$ and
$\dot{c}(0)$ is $\mathbf{g}_{\dot{c}(0)}$-orthogonal to $N$.
Denote by $J$ none-trivial $N$-Jacobi field along $c$ such that $J(0)\neq 0$. In
$\bar{M}^n$ we consider the same construction whose elements will be denoted by bar. We also suppose that $\mathbf{g}_{\dot{c}(t)}(J(t),J(t))|_0=
\mathbf{\bar{g}}_{\dot{\bar{c}}(t)}(\bar{J}(t),\bar{J}(t))|_0$.

Let us do a "transplant" of the vector field $J(t)$ in the manifold
$\bar{M}^n$ along the unit speed geodesic $\bar{c}(t)$, $0 \leqslant t \leqslant t^*$. We will follow [2].
Choose $\mathbf{g}_{\dot{c}(0)}$-orthonormal basis in
$T_{c(0)}M^n$ with $E_n = \dot{c}(0)$. Extend it to the parallel frame along $c(t)$. We obtain the basis $E_i(t)$ of
$\mathbf{g}_{\dot{c}(t)}$-orthonormal fields along $c(t)$ with
$E_n(t)=\dot{c}(t)$. We do the same in
$\bar{M}^n$ and obtain the $\mathbf{\bar{g}}_{\dot{\bar{c}}(t)}$-orthonormal fields $F_i(t)$ along $\bar{c}(t)$ with $F_n(t)=\dot{\bar{c}}(t)$.

Now, express $J(t)$ as $ \varphi^i(t)E_i(t)$ and define a new field $\tilde{J}(t) = \varphi^i(t)F_i(t)$. It satisfies the following properties [2, p.240]:
\begin{enumerate}
\item $J(t)$ and $\tilde{J}(t)$ have discontinuities at the same $t$ values;
\item
$\mathbf{g}_{\dot{c}(t)}(J(t),J(t))=\mathbf{\bar{g}}_{\dot{\bar{c}}(t)}(\tilde{J}(t),\tilde{J}(t))$;
\item
$\mathbf{g}_{\dot{c}(t)}(J(t),\dot{c}(t))=\mathbf{\bar{g}}_{\dot{\bar{c}}(t)}(\tilde{J}(t),\dot{\bar{c}}(t))$;
\item
$\mathbf{g}_{\dot{c}(t)}(J'(t),J'(t))=\mathbf{\bar{g}}_{\dot{\bar{c}}(t)}(\tilde{J}'(t),\tilde{J}'(t))$;
\item If $J(t^*)\neq 0$ and
$\mathbf{g}_{\dot{c}(t^*)}(J(t^*),\dot{c}(t^*))=0$, then
$\tilde{J}(t)$ can be chosen as $\tilde{J}(t^*) = \xi
\sqrt{\mathbf{g}_{\dot{c}(t^*)}(J(t^*),J(t^*))}$, where $\xi$ is a unit vector that is
$\mathbf{\bar{g}}_{\dot{\bar{c}}(t^*)}$-orthogonal to
$\dot{\bar{c}}(t^*)$.

\end{enumerate}

\textbf{Theorem 2.} \textit{Assume that for each $t\in
[0,s]$, for each flags $P\subset T_{c(t)}M^n$ and $\bar{P}\subset T_{\bar{c}(t)}\bar{M^n}$ along the geodesics
$c$, $\bar{c}$ respectively such that the flag $\bar{P}$ is the flag $P$ transplanted in $\bar{M^n}$ the inequality $K(\dot{c}(t),P)\leqslant
K(\dot{\bar{c}}(t),\bar{P})$ holds; the greater eigenvalue of the operator
$\overline{S}$ is less than or equal to the minimum eigenvalue of the operator $S$. Assume that there are no focal point along
 $\overline{c}$ to $\overline{N}$. Then for any $t\in
[0,s]$:
$$\mathbf{g}_{\dot{c}(t)}(J(t),J(t))\geqslant
\mathbf{\bar{g}}_{\dot{\bar{c}}(t)}(\bar{J}(t),\bar{J}(t))$$ and there are no focal points on
$c$. }

\textit{Proof.}
We give the sketch of proof. The idea repeats that of [2, p.245], [6], [12].

We prove that for each $t\in(0,s]$:

$$ \frac{\mathbf{g}_{\dot{c}(t)}(J'(t),J(t))}{\mathbf{g}_{\dot{c}(t)}(J(t),J(t))}\geqslant\frac{\mathbf{\bar{g}}_{\dot{\bar{c}}(t)}(\bar{J'}(t),\bar{J}(t))}{\mathbf{\bar{g}}_{\dot{\bar{c}}(t)}(\bar{J}(t),\bar{J}(t))} $$
 This is equivalent to the following

 $$\frac{d}{dt}\ln \frac{\mathbf{g}_{\dot{c}(t)}(J(t),J(t))}{\mathbf{\bar{g}}_{\dot{\bar{c}}(t)}(\bar{J}(t),\bar{J}(t))} \geqslant 0 $$
 and to the statement of the theorem.

Fix $t^*\in(0,s]$.

Choose  transplantation  $\tilde{J}$ such that $\tilde{J}(t^*) =
\frac{\sqrt{\bar{\textbf{g}}_{\dot{\bar{c}}(t)}(\bar{J}(t),\bar{J}(t))}}{\sqrt{\bar{\textbf{g}}_{\dot{\bar{c}}(t)}(\bar{J}(t),\bar{J}(t))}}|_{t^*} \bar{J}(t^*)$.

Denote $\lambda =
\mathbf{g}_{\dot{c}(t^*)}(J(t^*),J(t^*))$, $\bar{\lambda} =
\mathbf{\bar{g}}_{\dot{\bar{c}}(t^*)}(\bar{J}(t^*),\bar{J}(t^*))$.
Then

$$\frac{1}{\lambda}I_{t^*}(J,J)\geqslant
\frac{1}{\lambda}I_{t^*}(\tilde{J},\tilde{J})\geqslant
\frac{1}{\lambda}I_{t^*}(\frac{\sqrt{\lambda}}{\sqrt{\bar{\lambda}}}\bar{J},\frac{\sqrt{\lambda}}{\sqrt{\bar{\lambda}}} \bar{J}) = \frac{1}{\bar{\lambda}}I_{t^*}(\bar{J},\bar{J}).$$

The first inequality holds by virtue of curvature restrictions, conditions on the shape operators and of the properties of the transplanted field (see lemma 9.5.1, [2], p.242), the second one is just the Index Lemma (lemma 1).

From the Jacobi equation we get $I_{t^*}(J,J)
=\mathbf{g}_{\dot{c}(t^*)}(J(t^*),J'(t^*)) $. Since
$t^*$ is arbitrary the theorem follows. $\blacksquare$

\textbf{Corollary 1.} \textit{In Finsler manifold of flag curvature $K \leqslant -k^2$ and $\mathbf{T}$-curvature
$|\mathbf{T}| \leqslant \delta 0$, the exponential map of the hypersurface with all the normal curvatures
$\mathbf{k}_{\mathbf{n}}\geqslant \delta$ is non-degenerated and all the external parallel hypersurface are
regular.}

\section{On convexity of parallel hypersurfaces}

It is known that in Riemannian non-positive curved manifolds external parallel hypersurfaces to a locally convex hypersurface are locally convex. We prove analogous result with stronger restrictions on curvature.

\textbf{Lemma 2.} \textit{Assume that for a real $\lambda \geqslant 0$ a smooth function $f$ satisfies the following differential inequality}
$$f'(t)\leqslant f^2(t) - \lambda^2$$
\textit{Suppose that $f(0)\leqslant -\lambda$. Then  $f(t)\leqslant
-\lambda$ for all $t \geqslant 0$.}

\textit{Proof.} We prove that  $f$ cannot take value greater than  $-\lambda$. Otherwise, let $t_0>0$ be a point where $f(t_0)>-\lambda$. Denote by
$0<t_1<t_0$ the last point where $f(t_1)=-\lambda$. Then there exist the value $t_1<t_2<t_0$ such that $f'(t_2)$ would be greater than 0, $f(t_2)>-\lambda$ and consequently would not satisfy our differential inequality. And the lemma follows. $\square$

\textbf{Remark 1. } R. Bishop in [3] proved more stronger result in Riemannian geometry than proposition 1. He showed that if the second form is positively semidefinite then the hypersurface is locally convex. Thus the problem  \textit{does Bishop's result hold in general Finsler spaces} remains. Nevertheless for Berwald spaces Bishop's technique works using the affine property for Chern connection and in proposition 1 we can replace  strict inequalities by unstrict inequalities.

\textbf{Theorem 3.} \textit{ In Finsler manifold with flag curvature $K \leqslant -k^2$ and $\mathbf{T}$-curvature
$|\mathbf{T}| \leqslant \delta$ such that  $0\leqslant\delta < k$, consider a hypersurface $N$ with all the normal
curvatures $\mathbf{k}_{\mathbf{n}}> 2\delta$. Then all the external equidistant hypersurfaces are locally convex.}

\textit{Proof.} Denote by $N_t$ the external equidistant hypersurface to $N$, i.e.,
$N_t = \exp_N(t\nu(N))$. Such a mapping sends a point $x \in N$
to the point $\exp_x(t\mathbf{n})$. Consider the induced Riemannian metric $\tilde{g}$ associated with the geodesic vector field normal to $N$ at the point $x$. Proposition 2 implies equality of sectional curvature of induced metric and of the flag curvature of the initial Finsler metric for flags spanned on the normal directions. Ricatty equation on the normal curvatures of $N_t$ at corresponding points in the induced metric has the form
$$-\dot{\tilde{\mathbf{k}}}_{\mathbf{n}}(t) = \tilde{\mathbf{k}}_{\mathbf{n}}^2(t) +
g(t)$$

Here $g(t) \leqslant -k^2  <  -\delta^2$.

Proposition 3 yields
$\tilde{\mathbf{k}}_{\mathbf{n}}(0) > \delta$. By lemma 2 $\tilde{\mathbf{k}}_{\mathbf{n}}(t)
> \delta$ for all $t\geqslant  0$. Hence,
$\mathbf{k}_{\mathbf{n}}(t) \geqslant
\tilde{\mathbf{k}}_{\mathbf{n}}(t)-\delta > 0$, $t\geqslant  0$. And by proposition 1 the theorem follows.
 $\blacksquare$

\textbf{Remark 2. } Using remark 1 we can replace  strict inequalities by unstrict inequalities in theorem 3 for Berwald spaces.

\section{Proof of Theorem 1}

 Let $F = f(N^n)$ be a locally convex immersion of an $n$-dimensional compact manifold
$N^n$ in an $(n+1)$-dimensional complete simply connected Finsler manifold
$M^{n+1}$ with the imposed curvature conditions.

Define the external equidistant $F$ which will be denoted by $F_r$. It is regular by corollary 1. By theorem 3, $F_r$ is a locally convex hypersurface. Let
$S$ be a sphere in $M^{n+1}$ of a radius sufficiently large such that the ball bounded by this sphere contains
 $F_r$ and $F$ with small $r$. Proposition 4 implies convexity of the ball.

Consider the mapping $$\nu : F_r \rightarrow S$$ as the intersection point of the geodesic ray in the normal outer direction with the sphere
$S$. By virtue of non-existence of focal points (corollary 1) and of convexity of the sphere, locally different point are mapped to different. Hence $\nu$ is a local homeomorphism and therefore a covering. Since the sphere is simply connected for $n \geqslant 2$, the equidistance surface $F_r$ is homeomorphic to  $S$ and $\nu$ is a homeomorphism. Outside the ball rays normal to $F_r$ cannot intersect. Otherwise denote by $p$ the intersection point and consider the sphere $\tilde{S}$ in $M^{n+1}$ passing through the $p$ such that the ball bounded by this sphere contains
 $F_r$ and $F$ with small $r$ . Then we can construct the mapping $\tilde{\nu} : F_r \rightarrow \tilde{S}$. Since the mapping $\nu$ is a homeomorphism, the mapping $\tilde{\nu}$ is a homeomorphism too and we obtain the contrary.

Take the surface $F_{r+t_0}$ parallel to  $F_r$ and lying outside the ball bounded by $S$. Such a surface exists because of completeness of the ambient space. Locally it is an equidistant surface of $F$ hence it is locally convex. In addition it is an embedded surface without self-intersections, and it is the boundary of some body.

Arguing as in the proof of the Schmidt theorem [9], we can show that the body bounded by the surface  $F_{r+t_0}$ is convex and  $F_{r+t_0}$ is a compact hypersurface homeomorphic to the sphere (Schmidt theorem says that a connected set in $E^n$ having a local support flat at every boundary point is convex).

Now, consider the surfaces $F_{r+t_0-t}$ parallel to
$F_{r+t_0}$ and lying inside the body bounded by $F_{r+t_0}$. At small $t$ the surface $F_{r+t_0-t}$ is embedded, locally convex, and hence globally convex. Denote $t_1$ by the least upper bound of those
 $t$ for which the surfaces $F_{r+t_0-t}$ are embedded. Suppose that $t_1 < r+t_0$. The body bounded by
$F_{r+t_0-t_1}$ is convex, and  $F_{r+t_0-t_1}$
is a convex  surface. It must have points of self-tangency. The self-tangency may be internal or external. The former cannot occur because then the mapping $\nu$ is not a homeomorphism. The latter is also impossible, because the body bounded by $F_{r+t_0-t_1}$ is convex. This implies that $t_1 = r+t_0$ and $F$ is the boundary of convex  body. Considering a point inside this body we can show that hypersurface $F$ is homeomorphic to the sphere. $\blacksquare$

Using remark 2, one can prove Theorem 4.

\textbf{Theorem 4. }\textit{Let $f:N^n\rightarrow M^{n+1}$,
$n\geqslant 2$ be an immersion of compact connected manifold
$N^n$ in a complete simply connected Finsler manifold $M^{n+1}$ with Berwald metric of non-positive sectional curvature. If the immersion $f$ is locally  convex, then $f$ is an embedding, $f(N^n)$ is the boundary of a convex set in  $M^{n+1}$, and $f(N^n)$ is homeomorphic to the sphere
$\mathbb{S}^n$.}

\section{Example}

In [13] the example of non-complete simply connected negatively curved riemannian space which contains closed geodesic was constructed. If we consider the tube around the closed geodesic then we obtain a locally convex hypersurface (because the distance function from a convex set in non-positive curved riemannian manifold is convex), but the obtained hypersurface is evidently homeomorphic to the torus $S^n \times S^1$, here $(n+1)$ in a dimension of the ambient space. Thus the completeness condition is essential.

{99}
\end{document}